\documentclass{amsart}
\usepackage{amssymb,latexsym}
\topskip  \numberwithin{equation}{section}
\newtheorem{theorem}{Theorem}[section]

\newtheorem{corollary}[theorem]{Corollary}

\newtheorem{lemma}[theorem]{Lemma}

\begin{document}

\pagenumbering{arabic}
\pagestyle{headings}
\def\Des{\mbox{Des}}
\def\des{\mbox{des}}
\def\Ris{\mbox{Ris}}
\def\ris{\mbox{ris}}
\def\Lev{\mbox{Lev}}
\def\lev{\mbox{lev}}

\title{Counting descents, rises, and levels, with prescribed first element, in words}
\maketitle

\begin{center}
{\bf Sergey Kitaev}\\
{\it Institute of Mathematics, Reykjavik University,
IS-103 Reykjavik, Iceland}\\
{\tt sergey@ru.is}\\
\vskip 10pt
{\bf Toufik Mansour}\\
{\it Department of Mathematics, Haifa University, 31905 Haifa, Israel}\\
{\tt toufik@math.haifa.ac.il}\\
\vskip 10pt
{\bf Jeffrey B. Remmel}\\
{\it Department of Mathematics, University of California,
San Diego, La Jolla, CA 92093, USA}\\
{\tt jremmel@ucsd.edu}\\
\end{center}
\section*{Abstract}
Recently, Kitaev and Remmel~\cite{KR} refined the well-known
permutation statistic ``descent'' by fixing parity of one of the
descent's numbers. Results in~\cite{KR} were extended and
generalized in several ways in~\cite{HR, KR1, Liese, LR}. In this
paper, we shall fix a set partition of the natural numbers $\mathbb{N}$,
$(\mathbb{N}_1, \ldots, \mathbb{N}_t)$, and we study the distribution of descents, levels, and rises according to whether the first letter of the descent, rise, or level lies in $\mathbb{N}_i$ over the set
of words over the alphabet $[k]= \{1,\ldots,k\}$. In particular, we refine and
generalize some of the results in~\cite{BM}.
\section{Introduction}

The descent set of a permutation $\pi=\pi_1\ldots\pi_n\in S_n$ is
the set of indices $i$ for which $\pi_i>\pi_{i+1}$. This statistic
was first studied by MacMahon \cite{Mac} almost a hundred years ago
and it still plays an important role in the field of permutation
statistics. The number of permutations of length $n$ with exactly
$m$ descents is counted by the \emph{Eulerian number} $A_m(n)$. The
Eulerian numbers are the coefficients of the \emph{Eulerian
polynomials} $A_n(t)=\sum_{\pi\in S_n}t^{1+\des(\pi)}$. It is
well-known that the Eulerian polynomials satisfy the identity
$\sum_{m\geq0}m^nt^m=\frac{A_n(t)}{(1-t)^{n+1}}$. For more
properties of the Eulerian polynomials see \cite{Com}.

Recently, Kitaev and Remmel \cite{KR} studied the distribution of
a refined ``descent'' statistic on the set of permutations by fixing parity of
(exactly) one of the descent's numbers. For example, they showed
that the number of permutations in $S_{2n}$ (resp. $S_{2n+1}$) with
exactly $k$ descents such that the first entry of the descent is an
even number is given by $\binom{n}{k}^2n!^2$ (resp.
$\frac{1}{k+1}\binom{n}{k}^2(n+1)!^2$). In~\cite{KR1}, the authors
generalized results of~\cite{KR} by studying descents according to
whether the first or the second element in a descent pair is
equivalent to 0 mod $k\geq 2$.

Consequently, Hall and Remmel~\cite{HR} generalized results
of~\cite{KR1} by considering ``$X,Y$-descents,'' which are descents
whose ``top'' (first element) is in $X$ and whose ``bottom'' (second element)
is in $Y$  where $X$ and $Y$ are any sets of
the natural number $\mathbb{N}$. In particular, Hall and Remmel \cite{HR}
showed that one can reduce the problem of counting the number of permutations
$\sigma$ with $k$ $X,Y$-descents to the problem of computing
the $k$-th hit number of a Ferrers board in many cases.
Liese~\cite{Liese} also considered the situation of fixing equivalence
classes of both descent numbers simultaneously. Also,
papers~\cite{HLR} and~\cite{LR} discuss $q$-analogues of some of the
results in~\cite{HR, KR, KR1, Liese}. \\

Hall and Remmel \cite{HR} extended their results on counting
permutations with a given number of $X,Y$-descents to words. That
is, let $R(\rho)$ be the rearrangement class of the word
$1^{\rho_1}2^{\rho_2}\cdots m^{\rho_m}$ (i.e., $\rho_1$ copies of
$1$, $\rho_2$ copies of $2$, etc.) where $\rho_1 + \cdots + \rho_m =n$. For any set $X \subseteq
\mathbb{N}$ and any set $[m]= \{1, 2,\ldots, m\}$, we let $X_m  = X
\cap [m]$ and $X_m^c = [m]-X$. Then given $X, Y \subseteq
\mathbb{N}$  and a word $w = w_1 \cdots w_n \in R(\rho)$, define
\begin{eqnarray*}
Des_{X,Y}(w) & = & \{ i : w_i > w_{i+1}
~\&~ w_i \in X ~\&~ w_{i+1} \in Y \} ,\\
des_{X,Y}(w) & = & |Des_{X,Y}(w)|, \mbox{ and} \\
P_{\rho,s}^{X,Y} & = & \left| \{ w \in R(\rho) :des_{X,Y}(w) = s \} \right|.
\end{eqnarray*}
Hall and Remmel \cite{HR} proved the following theorem by purely combinatorial
means.
\begin{theorem}\label{combword1}
\begin{equation}\label{eq:combword1}
P_{\rho,s}^{X,Y} = {a \choose \rho_{v_1},\rho_{v_2}, \ldots, \rho_{v_b}} \sum\limits_{r=0}^s (-1)^{s-r} {a+r \choose r}
{n+1 \choose s-r}
\prod\limits_{x\in X} {\rho_{x} + r + \alpha_{X,\rho,x}+\beta_{Y,\rho,x}  \choose \rho_{x} },
\end{equation}
where $X_m^c = \{v_1, v_2, \ldots, v_b \}, a = \sum\limits_{i=1}^b \rho_{v_i}$, and for any $x \in X_m$,
\begin{eqnarray*}
\alpha_{X,\rho,x} & = & \sum\limits_{\tiny \begin{array}{c} z \notin X \\ x < z \leq m \end{array} \normalsize}
\rho_z, \mbox{ and} \\
\beta_{Y,\rho,x} & = & \sum\limits_{\tiny \begin{array}{c} z \notin Y \\ 1 \leq z < x \end{array} \normalsize} \rho_z.
\end{eqnarray*}
\end{theorem}

In this paper, we shall study similar statistics over the set
$[k]^n$ of $n$-letter words over fixed finite alphabet $[k] =\{1,
2,\ldots, k\}$. In what follows, $E=\{2,4,6,\ldots\}$ and
$O=\{1,3,5,\ldots\}$ are the sets of even and odd numbers
respectively. Also, we let ${\bf x}[t]=(x_1,\ldots,x_t)$. Then given
a word $\pi=\pi_1\pi_2\ldots\pi_n\in [k]^n$ and a set $X \subseteq
\mathbb{N}$, we define the following statistics:
\begin{itemize}
\item $\overleftarrow{\Des}_X(\pi)=\{i: \pi_i>\pi_{i+1}\mbox{ and }\pi_i\in
X\}$ and
$\overleftarrow{\des}_X(\pi)=|\overleftarrow{\Des}_X(\pi)|$,\\[-7pt]

\item $\overleftarrow{\Ris}_X(\pi)=\{i: \pi_i<\pi_{i+1}\mbox{ and }\pi_{i}\in
X\}$ and $\overleftarrow{\ris}_X(\pi)=|\overleftarrow{\Ris}_X(\pi)|$,\\[-7pt]

\item $\Lev_X(\pi)=\{i: \pi_i=\pi_{i+1}\mbox{ and }\pi_{i}\in
X\}$ and $\lev_X(\pi)=|\Lev_X(\pi)|$.\\[-7pt]
\end{itemize}
Let $(\mathbb{N}_1, \ldots, \mathbb{N}_t)$ be a set partition of the
natural numbers $\mathbb{N}$, i.e. $\mathbb{N}=\mathbb{N}_1\cup
\mathbb{N}_2\cup\ldots \cup \mathbb{N}_t$ and $\mathbb{N}_i\cap
\mathbb{N}_j=\emptyset$ for $i\neq j$. Then the main goal of this
paper is to study the following multivariate generating function
(MGF)
\begin{equation}\label{eq:MGF}
A_k=A_k({\bf x}[t];{\bf y}[t];{\bf z}[t];{\bf q}[t])=\\
\sum_{\pi}\prod_{i=1}^{t}x_i^{\overleftarrow{\des}_{\mathbb{N}_i}(\pi)}
y_i^{\overleftarrow{\ris}_{\mathbb{N}_i}(\pi)}z_i^{\lev_{\mathbb{N}_i}(\pi)}q_i^{i(\pi)}
\end{equation}
where $i(\pi)$ is the number of letters from $\mathbb{N}_i$ in $\pi$
and the sum is over all words over $[k]$.

The outline of this paper is as follows. In section~\ref{sec2}, we
shall develop some general methods to compute (\ref{eq:MGF}). In
section~\ref{sec3}, we shall concentrate on the computing generating
functions for the distribution of the number of levels. That is, we
shall study $A_k$ where set $x_i=y_i=1$ for all $i$. In
section~\ref{sec4}, we shall find formulas for the number of words
in $[k]^n$ that have $s$ descents that start with an element less
than or equal to $t$ (greater than $t$) for any $t \leq k$. Note
that if we replace a word $w = w_1 \cdots w_n \in [k]^n$ by its
complement $w^c = (k+1 -w_1) \cdots (k+1 - w_n)$, then it is easy to
see that $\overleftarrow{\des}_{[t]}(w)  =
\overleftarrow{\ris}_{\{k+1-t, \ldots, k\}}(w)$ and
$\overleftarrow{\des}_{\{t+1, \ldots, k\}}(w)  =
\overleftarrow{\ris}_{[k-t]}(w)$. Thus we will also obtain formulas
for the the number of words in $[k]^n$ that have $s$ rises that
start with an element less than or equal to $t$ (greater than $t$)
for any $t \leq k$. We will also show that in the cases where $t=2$
or $t=k-1$, there are alternative ways to compute our formulas which
lead to non-trivial binomial identities. In section~\ref{sec5}, we
shall apply our results to study the problem of counting the number
of words in $[k]^n$ with $p$ descents (rises) that start with an
element which is equivalent to $i \mod s$ for any $s \geq 2$ and $i
=1, \ldots, s$. In particular, if $s \geq 2$ and $(\mathbb{N}_1,
\ldots, \mathbb{N}_s)$ is the set partition of $\mathbb{N}$ where
$\mathbb{N}_i = \{x \in \mathbb{N}: x \equiv i \mod s\}$ for $i =1,
\ldots, s$, then we shall study the generating functions
\begin{equation}\label{ask}
A^{(s)}_k({\bf x}[s];{\bf y}[s];{\bf z}[s];{\bf q}[s]) =
\sum_{\pi}\prod_{i=1}^{s}x_i^{\overleftarrow{\des}_{\mathbb{N}_i}(\pi)}
y_i^{\overleftarrow{\ris}_{\mathbb{N}_i}(\pi)}z_i^{\lev_{\mathbb{N}_i}(\pi)}q_i^{i(\pi)}
\end{equation}
and
\begin{equation}\label{askq}
A^{(s)}_k({\bf x}[s];{\bf y}[s];{\bf z}[s];q) = \sum_{\pi} q^{|\pi|}
\prod_{i=1}^{s}x_i^{\overleftarrow{\des}_{\mathbb{N}_i}(\pi)}
y_i^{\overleftarrow{\ris}_{\mathbb{N}_i}(\pi)}z_i^{\lev_{\mathbb{N}_i}(\pi)}.
\end{equation}
Our general results in section 2 allow us to derive an explicit formula
for $A^{(s)}_k({\bf x}[s];{\bf y}[s];{\bf z}[s];{\bf q}[s])$ depending on
the equivalence class of $k$ mod $s$.

For example, in the case where $s =2$, our general result implies that
\begin{eqnarray}\label{eq:MFG1EO}
A^{(2)}_{2k}(q_1,q_2) &=& A_k(x_1,x_2,y_1,y_2,z_1,z_2,q_1,q_2)  = \\
&=& \sum_{\pi}
x_1^{\overleftarrow{\des}_O(\pi)}x_2^{\overleftarrow{\des}_E(\pi)}
y_1^{\overleftarrow{\ris}_O(\pi)}y_2^{\overleftarrow{\ris}_E(\pi)}
z_1^{\lev_O(\pi)}z_2^{\lev_E(\pi)}q_1^{\mbox{odd}(\pi)}q_2^{\mbox{even}(\pi)}
\nonumber \\
&=& \frac{1 + (\lambda_1 \mu_2 + \lambda_2) \frac{1-\mu_1^k\mu_2^k}{1-\mu_1\mu_2}}{1 - (\nu_1 \mu_2 + \nu_2) \frac{1-\mu_1^k\mu_2^k}{1-\mu_1\mu_2}}
\nonumber
\end{eqnarray}
where the sum is over all words over $[2k]$, even$(\pi)$ (resp.
odd$(\pi)$) is the number of even (resp. odd) numbers in $\pi$,
$\lambda_j = \frac{q_j(1-y_j)}{1-q_j(z_j-y_j)}$, $\mu_i =
\frac{q_i(z_i-x_i)}{1-q_i(z_i-y_i)}$, and $\nu_j =
\frac{q_jy_j}{1-q_j(z_j-y_j)}$ for $j=1,2$. Then by specializing the
variables appropriately, we will find explicit formulas for the
number of words $w \in [2k]^n$ such that
$\overleftarrow{\des}_E(\pi) = p$, $\overleftarrow{\des}_O(\pi) =
p$, $\overleftarrow{\ris}_E(\pi) = p$, $\overleftarrow{\ris}_O(\pi)
= p$, etc. For example, we prove that the number of $n$-letter words
$\pi$ on $[2k]$ having $\overleftarrow{\des}_O(\pi)=p$ {\rm(}resp.
$\overleftarrow{\ris}_E(\pi)=p${\rm)} is given by
$$\sum_{j=0}^n\sum_{i=0}^j(-1)^{n+p+i}2^j\binom{j}{i}\binom{i}{n}
\binom{n-j}{p}.$$ In fact, we shall show that similar formulas hold
for the number of words $\pi \in [k]^n$ with $p$ descents (rises,
levels) whose first element is equivalent to $t \mod s$ for any $s
\geq 2$ and $0 \leq t \leq s-1$. Our results refine and generalize
the results in~\cite{BM} related to the distribution of descents,
levels, and rises in words.  Finally, in section~\ref{sec6}, we
shall discuss some open questions and further research.

\section{The general case}\label{sec2}

We need the following notation:
$$A_k(i_1,\ldots,i_m)=A_k({\bf x}[t];{\bf y}[t];{\bf z}[t];{\bf q}[t];{\bf i}[m])=
\sum_{\pi}\prod_{i=1}^{t}x_i^{\overleftarrow{\des}_{\mathbb{N}_i}(\pi)}
y_i^{\overleftarrow{\ris}_{\mathbb{N}_i}(\pi)}z_i^{\lev_{\mathbb{N}_i}(\pi)}q_i^{i(\pi)}$$
where the sum is taken over all words $\pi=\pi_1\pi_2\ldots$ over
$[k]$ such that $\pi_1\ldots\pi_m=i_1\ldots i_m$.

From our definitions, we have that
\begin{equation}\label{eqbas}
A_k=1+\sum_{i=1}^kA_k(i).
\end{equation}
Thus, to find a formula for $A_k$, it is sufficient to find a formula
for $A_k(i)$ for each $i=1,2,\ldots,k$. First let us find a
recurrence relation for the generating function $A_k(i)$.

\begin{lemma}\label{lemg1}
For each $s\in \mathbb{N}_i$, $1\leq s\leq k$ and $1\leq i\leq t$,
we have
\begin{equation}\label{eq:lemg1}
A_k(s)=\frac{q_iy_i}{1-q_i(z_i-y_i)}A_k+\frac{q_i(1-y_i)}{1-q_i(z_i-y_i)}+\frac{q_i(x_i-y_i)}{1-q_i(z_i-y_i)}
\sum_{j=1}^{s-1}A_k(j).
\end{equation}
\end{lemma}
\begin{proof}
From the definitions we have that
$$\begin{array}{ll}
A_k(s)&=q_i+\sum_{j=1}^{k}A_k(s,j)\\[4pt]
&=q_i+\sum_{j=1}^{s-1}A_k(s,j)+A_k(s,s)+\sum_{j=s+1}^{k}A_k(s,j).
\end{array}$$
Let $\pi$ be any $n$-letter word over $[k]$ where $n\geq 2$ and
$\pi_1=s>\pi_2=j$. If we let $\pi'=\pi_2\pi_3\ldots\pi_n$, then it is
easy to see that
$$\begin{array}{ll}
\overleftarrow{\des}_{\mathbb{N}_i}(\pi)=1+\overleftarrow{\des}_{\mathbb{N}_i}(\pi'),&
i(\pi)=1+i(\pi').\end{array}$$ It is also easy to see
that remaining $4t-2$ statistics of interest take the same value on
$\pi$ and $\pi'$.

This implies that $A_k(s,j)=q_ix_iA_k(j)$ for each $1\leq j<s$.
Similarly, $A_k(s,s)=q_iz_iA_k(s)$ and $A_k(s,j)=q_iy_iA_k(j)$ for
$s< j\leq k$. Therefore,
$$A_k(s)=q_i+q_ix_i\sum_{j=1}^{s-1}A_k(j)+q_iz_iA_k(s)+q_iy_i\sum_{j=s+1}^{k}A_k(j).$$
Using \eqref{eqbas}, we have
$\sum_{j=s+1}^kA_k(j)=A_k-\sum_{j=1}^{s-1}A_k(j)-A_k(s)-1$, and thus
$$A_k(s)=\frac{q_iy_i}{1-q_i(z_i-y_i)}A_k+\frac{q_i(1-y_i)}{1-q_i(z_i-y_i)}+\frac{q_i(x_i-y_i)}{1-q_i(z_i-y_i)}
\sum_{j=1}^{s-1}A_k(j),$$ as desired.
\end{proof}

\begin{lemma}\label{lemg2}
For each $k \geq 1$ and $s \in [k]$,
\begin{equation}\label{eq:lemg2}
\sum_{j=1}^s A_k(j) = \sum_{j=1}^s \gamma_j \prod_{i=j+1}^s (1-\alpha_i)
\end{equation}
where, for $i\in\mathbb{N}_m$ and $i\geq 1$,
$\gamma_i=\frac{q_my_m}{1-q_m(z_m-y_m)}A_k+\frac{q_m(1-y_m)}{1-q_m(z_m-y_m)}$
and $\alpha_i=\frac{q_m(y_m-x_m)}{1-q_m(z_m-y_m)}$.
\end{lemma}
\begin{proof}
We proceed by induction on $s$. Note, that
given our definitions of $\gamma_i$ and $\alpha_i$, we can rewrite
(\ref{eq:lemg1}) as
\begin{equation}\label{eq:lemg22}
A_k(s) = \gamma_s -\alpha_s\sum_{j=1}^{s-1} A_k(j).
\end{equation}
It follows that
$$A_k(1) = \gamma_1$$
so that (\ref{eq:lemg2}) holds for $s =1$. Thus the base case of our induction
holds. Now assume that
(\ref{eq:lemg2}) holds for $s$ where $1 \leq s < k$. Then using our induction
hypothesis and (\ref{eq:lemg22}), it follows that
\begin{eqnarray*}
&&A_k(1) + \cdots + A_k(s) + A_k(s+1)\\
&&\quad =\sum_{j=1}^s \gamma_j \prod_{i=j+1}^s (1-\alpha_i) +
\gamma_{s+1} -
\alpha_{s+1}\left( \sum_{j=1}^s \gamma_j \prod_{i=j+1}^s (1-\alpha_i)\right)\\
&&\quad =\gamma_{s+1} - \sum_{j=1}^s \gamma_j \prod_{i=j+1}^{s+1} (1-\alpha_i)\\
&&\quad=\sum_{j=1}^{s+1} \gamma_j \prod_{i=j+1}^{s+1} (1-\alpha_i).
\end{eqnarray*}
 Thus the induction step also holds so that (\ref{eq:lemg2}) must hold in
general.
\end{proof}

Lemma~\ref{lemg1} gives that $A_k(i)$, for $1\leq i\leq k$, are the
solution to the following matrix equation
\begin{equation}\label{eqmain}\left(\begin{array}{lllllll}
1     &0     &0     &0&0&\ldots&0\\
\alpha_2&1     &0     &0&0&\ldots&0\\
\alpha_3 &\alpha_3 &1     &0&0&\ldots&0\\
\alpha_4&\alpha_4&\alpha_4&1&0&\ldots&0\\
\vdots&&&\vdots&&&\vdots\\
\alpha_k&\alpha_k&\alpha_k&\alpha_k&\alpha_k&\ldots&1
\end{array}\right)
\cdot\left(\begin{array}{l}A_{k}(1)\\A_{k}(2)\\\vdots\\A_{k}(k)\end{array}\right)
=\left(\begin{array}{l}\gamma_1\\\gamma_2\\\vdots\\\gamma_k\end{array}\right)\end{equation}
where, for $i\in\mathbb{N}_m$ and $i\geq 1$,
$\gamma_i=\frac{q_my_m}{1-q_m(z_m-y_m)}A_k+\frac{q_m(1-y_m)}{1-q_m(z_m-y_m)}$,
and, for $i\in\mathbb{N}_m$ and $i\geq 2$,
$\alpha_i=\frac{q_m(y_m-x_m)}{1-q_m(z_m-y_m)}$. Notice that
$\alpha_i=\alpha_j$ and $\gamma_i=\gamma_j$ whenever $i$ and $j$ are
from the same set $\mathbb{N}_m$ for some $m$.
In fact, it is easy to see that (\ref{eq:lemg1}) and (\ref{eq:lemg2})
imply that
\begin{equation}\label{eq:lemg3}
A_k(i)=\gamma_i-\alpha_i\sum_{j=1}^{i-1}\gamma_j\prod_{i=j+1}^{i-1}(1-\alpha_i)
\end{equation}
holds for $i=1, \ldots, k$ so that (\ref{eqmain}) has an explicit solution.
By combining (\ref{eqbas}) and (\ref{eq:lemg3}), we can obtain the following
result.

\begin{theorem}\label{main} For $\alpha_i$ and $\gamma_i$ as above {\rm(}defined right below
\eqref{eqmain}{\rm)}, we have
$$A_k=1+\sum_{j=1}^k\gamma_j\prod_{i=j+1}^k(1-\alpha_i)$$ solving which for $A_k$ gives
\begin{equation}\label{eq:Main}
A_k=\frac{1+\sum_{j=1}^{k}\frac{q_j(1-y_j)}{1-q_j(z_j-y_j)}\prod_{i=j+1}^{k}\frac{1-q_i(z_i-x_i)}{1-q_i(z_i-y_i)}}
{1-\sum_{j=1}^{k}\frac{q_jy_j}{1-q_j(z_j-y_j)}\prod_{i=j+1}^{k}\frac{1-q_i(z_i-x_i)}{1-q_i(z_i-y_i)}}
\end{equation}
where for each variable $a\in\{x,y,z,q\}$ we have $a_i=a_m$ if
$i\in\mathbb{N}_m$.
\end{theorem}

Even though we state Theorem~\ref{main} as the main theorem in this
paper, its statement can be (easily) generalized if one considers
compositions instead of words. Indeed, let

$$B_k=B_k({\bf x}[t];{\bf y}[t];{\bf z}[t];{\bf q}[t];v)=\\
\sum_{\pi}v^{|\pi|}\prod_{i=1}^{t}x_i^{\overleftarrow{\des}_{\mathbb{N}_i}(\pi)}
y_i^{\overleftarrow{\ris}_{\mathbb{N}_i}(\pi)}z_i^{\lev_{\mathbb{N}_i}(\pi)}q_i^{i(\pi)}$$
where the sum is taken over all compositions $\pi=\pi_1\pi_2\ldots$
with parts in $[k]$ and $|\pi|=\pi_1+\pi_2+\cdots$ is the weight of
the composition $\pi$. Also, we let
$$B_k(i_1,\ldots,i_m)=B_k({\bf
x}[t];{\bf y}[t];{\bf z}[t];{\bf q}[t];{\bf i}[m]; v)=
\sum_{\pi}v^{|\pi|}\prod_{i=1}^{t}
x_i^{\overleftarrow{\des}_{\mathbb{N}_i}(\pi)}
y_i^{\overleftarrow{\ris}_{\mathbb{N}_i}(\pi)}
z_i^{\lev_{\mathbb{N}_i}(\pi)}q_i^{i(\pi)}$$
where again the sum is taken over all compositions
$\pi=\pi_1\pi_2\ldots$ with parts in $[k]$.

Next, one can copy the arguments of Lemma~\ref{lemg1} substituting
$q_i$ by $v^sq_i$ to obtain the following generalization of Lemma~\ref{lemg1}:
$$B_k(s)=\frac{v^sq_iy_i}{1-q_i(z_i-y_i)}B_k+\frac{v^sq_i(1-y_i)}{1-q_i(z_i-y_i)}+\frac{v^sq_i(x_i-y_i)}{1-q_i(z_i-y_i)}
\sum_{j=1}^{s-1}B_k(j).$$ One can then prove the obvious analogue of
Lemma~\ref{lemg1} by induction and apply it to prove the following theorem.

\begin{theorem}\label{main1} We have
$$B_k=1+\sum_{j=1}^k\gamma_j\prod_{i=j+1}^k(1-\alpha_i)$$ where
$\gamma_i=\frac{v^iq_my_m}{1-q_m(z_m-y_m)}B_k+\frac{v^iq_m(1-y_m)}{1-q_m(z_m-y_m)}$,
and $\alpha_i=\frac{v^iq_m(y_m-x_m)}{1-q_m(z_m-y_m)}$ if $i$ belongs
to $\mathbb{N}_m$. Thus,
$$B_k=\frac{1+\sum_{j=1}^{k}\frac{v^jq_j(1-y_j)}{1-q_j(z_j-y_j)}\prod_{i=j+1}^{k}\frac{1-q_i(z_i-y_i+v^i(y_i-x_i))}{1-q_i(z_i-y_i)}}
{1-\sum_{j=1}^{k}\frac{v^jq_jy_j}{1-q_j(z_j-y_j)}\prod_{i=j+1}^{k}\frac{1-q_i(z_i-y_i+v^i(y_i-x_i))}{1-q_i(z_i-y_i)}}$$
where for each variable $a\in\{x,y,z,q\}$ we have $a_i=a_m$ if $i\in
\mathbb{N}_m$.
\end{theorem}
Theorem~\ref{main1} can be viewed as a $q$-analogue to
Theorem~\ref{main}. (Set $v=1$ in Theorem~\ref{main1} to get
Theorem~\ref{main}.)

\section{Counting words by the types of levels}\label{sec3}

Suppose we are given a set partition
$\mathbb{N}=\mathbb{N}_1\cup\mathbb{N}_2 \cup \cdots \cup
\mathbb{N}_s$.  First observe that for any fixed $i$, if we want the
distribution of words in $[k]^n$ according to the number of levels
which involve elements in $\mathbb{N}_i$, then it is easy to see by
symmetry that the distribution  will depend only on the cardinality
of $\mathbb{N}_i \cap [k]$. Thus we only need to consider the case
where $s=2$ and $\mathbb{N}_1 = \{1,\ldots,t\}$ for some $t \leq k$.

Let
\begin{eqnarray}
\lambda_j &=& \frac{q_j(1-y_j)}{1-q_j(z_j-y_j)}, \\
\nu_j &=& \frac{q_jy_j}{1-q_j(z_j-y_j)}, \ \mbox{and}\\
\mu_i &=& \frac{1-q_i(z_i-x_i)}{1-q_i(z_i-y_i)}.
\end{eqnarray}
Then we can rewrite (\ref{eq:Main}) for any arbitrary set partition
$\mathbb{N}=\mathbb{N}_1\cup\mathbb{N}_2 \cup \cdots \cup
\mathbb{N}_s$ as
\begin{equation}\label{eq:3.1}
A_k= \frac{1+ \sum_{j=1}^k \lambda_j \prod_{i=j+1}^k \mu_i}{1- \sum_{j=1}^k \nu_j \prod_{i=j+1}^k \mu_i}
\end{equation}
where for each variable $a \in \{x,y,z,q\}$, we have $a_i = a_m$ if
$i \in \mathbb{N}_m$.

Suppose we set $x_1 = x_2 = y_1 =y_2 = z_2 =1$ and $q_1 =q_2 = q$ in
(\ref{eq:3.1}) in the special case where $s=2$ and $\mathbb{N}_1
=[t]$ for some $t \leq k$. Then $\lambda_1 = \lambda_2 =0$, $\nu_1 =
\frac{q}{1-q(z_1-1)}$, $\nu_2 =q$, and $\mu_1 = \mu_2 =1$. It
follows that in this case,
\begin{eqnarray*}
A_k &=& \frac{1}{1-\left(\frac{tq}{1-q(z_1-1)} + q(k-t)\right)}\\
&=& \sum_{m \geq 0} q^m \left(\frac{t}{1-q(z_1-1)} + (k-t)\right)^m \\
&=& \sum_{m \geq 0} q^m \sum_{i=0}^m \binom{m}{i} (k-t)^{m-i} t^i
\left(\frac{1}{1-q(z_1-1)}\right)^i.
\end{eqnarray*}
Since
\begin{eqnarray}\label{eq:3.2}
\left(\frac{1}{1-q(z_1-1)}\right)^i &=& \sum_{a \geq 0} \frac{(i)_a}{a!}q^a
(z_1-1)^a \nonumber \\
&=& \sum_{a \geq 0} \binom{i+a-1}{a}q^a (z_1-1)^a,
\end{eqnarray}
it follows that
\begin{equation}\label{eq:3.3}
A_k = \sum_{n \geq 0} q^n \sum_{m=0}^n \sum_{i=0}^m
\binom{m}{i} \binom{i+n-m-1}{n-m} (k-t)^{m-i} t^i (z_1-1)^{n-m}.
\end{equation}
Thus taking the coefficient of $z_1^s$ on both sides of
(\ref{eq:3.3}), we obtain the following result.

\begin{theorem}\label{thm:3.1}
Let $\mathbb{N}=\mathbb{N}_1\cup\mathbb{N}_2$ where $\mathbb{N}_1 =[t]$ and
$\mathbb{N}_2 = \mathbb{N} -\mathbb{N}_1$. Then if $t \leq k$, the number of
words in $[k]^n$ with $s$ levels that start with elements in
$\mathbb{N}_1$ is
\begin{equation}\label{eq:3.4}
\sum_{m=0}^n \sum_{i=0}^m (-1)^{n-m-s}
\binom{m}{i} \binom{i+n-m-1}{n-m} \binom{n-m}{s} (k-t)^{m-i} t^i.
\end{equation}
\end{theorem}

Going back to the general set partition
$\mathbb{N}=\mathbb{N}_1\cup\mathbb{N}_2 \cup \cdots \cup
\mathbb{N}_s$, we can obtain a general formula for the number of
words in $[k]^n$ for which there are $t_i$ levels which start with
an element of $\mathbb{N}_i$ for $i =1,\ldots, s$ as follows. Let
$n_i = |\mathbb{N}_i \cap [k]|$ for $i =1, \ldots, n$. Then if set
$x_j =y_j =1$ and $q_j =q$ for all $j$, then it will be the case
that $\lambda_j =0$ and $\mu_j =1$ and $\nu_j = \frac{q}{1-q(z_j
-1)}$ for all $j$. It easy follows that in this case,
\begin{eqnarray*}
A_k &=& \frac{1}{1-\left( \sum_{i=1}^s \frac{n_iq}{1-q(z_i -1)}\right)} \\
&=& \sum_{m \geq 0}q^m \left( \sum_{i=1}^s \frac{n_i}{1-q(z_i -1)}\right)^m \\
&=& \sum_{m \geq 0}q^m \sum_{\overset{a_1+ \cdots a_s =m}{a_1, \ldots, a_s \geq 0}} \binom{m}{a_1,\ldots, a_m} \prod_{i=1}^s \left( \frac{n_i}{1-q(z_i-1)}\right)^{a_i}.
\end{eqnarray*}
Then using (\ref{eq:3.2}), we see that
\begin{eqnarray}\label{eq:3.5}
A_k &=& \sum_{m \geq 0}q^m \sum_{\overset{a_1+ \cdots +a_s =m}{a_1,
\ldots, a_s \geq 0}} \binom{m}{a_1,\ldots, a_m} n_1^{a_1} \cdots
n_s^{a_s} \prod_{i=1}^s \sum_{b_i \geq 0}
\frac{(a_i)_{b_i}}{b_i!} q^{b_i} (z_1 -1)^{b_i} \nonumber \\
&=& \sum_{n \geq 0} q^n \sum_{m=0}^n \sum_{\overset{a_1+ \cdots +a_s
=m}{a_1, \ldots, a_s \geq 0}} \sum_{\overset{b_1+ \cdots +b_s
=n-m}{b_1, \ldots, b_s \geq 0}} \binom{m}{a_1,\ldots, a_m} n_1^{a_1}
\cdots n_s^{a_s} \prod_{i=1}^s \binom{a_i + b_i -1}{b_i} (z_i
-1)^{b_i}.
\end{eqnarray}
Taking the coefficient of $z_1^{t_1} \cdots z_s^{t_s}$ on both sides of
(\ref{eq:3.5}), we obtain the following result.

\begin{theorem}\label{thm:3.2}
Let $\mathbb{N}=\mathbb{N}_1\cup \cdots \cup \mathbb{N}_s$ be a set
partition of $\mathbb{N}$. Let $n_i = |\mathbb{N}_i \cap [k]|$ for
$i =1, \ldots, s$. Then the number of words in $[k]^n$ with $t_i$
levels that start with elements in $\mathbb{N}_i$ for $i =1, \ldots,
s$ is
\begin{equation}\label{eq:3.6}
\sum_{m=0}^n \sum_{\overset{a_1+ \cdots +a_s =m}{a_1, \ldots, a_s
\geq 0}} \sum_{\overset{b_1+ \cdots +b_s =n-m}{b_1, \ldots, b_s \geq
0}} \binom{m}{a_1,\ldots, a_m} n_1^{a_1} \cdots n_s^{a_s}
\prod_{i=1}^s \binom{a_i + b_i -1}{b_i}\binom{b_i}{t_i}.
\end{equation}
\end{theorem}

\section{Classifying words by the number of descents that start with elements $\leq t$ ($\geq t+1$).}\label{sec4}

In this section, we shall consider the set partition $\mathbb{N} = \mathbb{N}_1 \cup \mathbb{N}_2$ where $\mathbb{N}_1 =[t]$. Now if $t \leq k$, then it
is easy to see that we can rewrite (\ref{eq:Main}) as
\begin{equation}\label{eq:4.1}
A_k= \frac{1+ \sum_{j=1}^k \lambda_j \prod_{i=j+1}^k \mu_i}{1- \sum_{j=1}^k \nu_j \prod_{i=j+1}^k \mu_i}
\end{equation}
where
$$
\lambda_j =\frac{q_1(1-y_1)}{1-q_1(z_1-y_1)},\quad \nu_j = \frac{q_1
y_1}{1-q_1(z_1-y_1)},\mbox{ and }\mu_j =
\frac{q_1(z_1-x_1)}{1-q_1(z_1-y_1)}\mbox{ if $j \leq t$}$$
and
$$\lambda_j = \frac{q_2(1-y_2)}{1-q_2(z_2-y_2)},\quad\nu_j = \frac{q_2
y_2}{1-q_2(z_2-y_2)},\mbox{ and }\mu_j =
\frac{q_2(z_2-x_2)}{1-q_2(z_2-y_2)}\mbox{ if $j > t$}.$$ Now if we
want to find formulas for the number of words in $[k]^n$ with $s$
descents that start with an element less than or equal to $t$, then
we need to set $x_2 = y_1 = y_2 = z_1 = z_2 =1$ and $q_1 =q_2 =q$ in
(\ref{eq:4.1}). In that case, we will have $\lambda_j = 0$  and
$\nu_j = q$ for all $j$, $\mu_j = 1+q(x_1-1)$ for $j \leq t$, and
$\mu_j = 1$ for $j > t$. It follows that
\begin{eqnarray*}
A_k &=& \frac{1}{1-\left( \sum_{j=t+1}^k q + \sum_{j=1}^t q\prod_{i=j+1}^t (1+q(x_1 -1)) \right)} \\
&=& \frac{1}{1-\left( (k-t)q + q\frac{(1+q(x_1-1))^t -1}{(1+q(x_1-1)) -1}\right)}\\
&=& \frac{1}{1-\frac{1}{(x_1-1)}\left( (k-t)q(x_1-1)-1 + (1+q(x_1-1))^t \right)}\\
&=& \sum_{m \geq 0} \frac{1}{(x_1-1)^m} \left( (k-t)q(x_1-1)-1 + (1+q(x_1-1))^t \right)^m \\
&=& \sum_{m \geq 0} \frac{1}{(x_1-1)^m}\sum_{a=0}^m \binom{m}{a}
((k-t)q(x_1-1)-1)^{m-a} (1+q(x_1-1))^{ta}\\
&=& \sum_{m \geq 0} \frac{1}{(x_1-1)^m}\sum_{a=0}^m \sum_{b=0}^{m-a}
\sum_{c=0}^{ta} \binom{m}{a} \binom{m-a}{b} \binom{ta}{c}
(-1)^{m-a-b} q^b(k-t)^b(x_1-1)^b q^c (x_1-1)^c.
\end{eqnarray*}
If we want to take the coefficient of $q^n$, then we must have
$b+c =n$ or $c = n-b$. Thus
\begin{equation}\label{eq:4.2}
A_k = \sum_{n \geq 0} q^n \sum_{m \geq 0} \sum_{a=0}^m \sum_{b=0}^{m-a} \binom{m}{a} \binom{m-a}{b} \binom{ta}{n-b}
(-1)^{m-a-b} (k-t)^b(x_1-1)^{n-m}.
\end{equation}
Taking the coefficient of $q^n$ of both sides of
(\ref{eq:4.2}), we see that
\begin{equation}\label{eq:4.3}
\sum_{\pi \in [k]^n} x_1^{\overleftarrow{\des}_{[t]}(\pi)} = \sum_{m
\geq 0} \sum_{a=0}^m \sum_{b=0}^{m-a} \binom{m}{a} \binom{m-a}{b}
\binom{ta}{n-b} (-1)^{m-a-b} (k-t)^b(x_1-1)^{n-m}
\end{equation}
for all $n$. However, if we replace $x_1$ by $z+1$ in (\ref{eq:4.3}), we
see that the polynomial
$$\sum_{\pi \in [k]^n} (z+1)^{\overleftarrow{\des}_{[t]}(\pi)}$$
has the Laurent expansion
$$
\sum_{m \geq 0} \sum_{a=0}^m \sum_{b=0}^{m-a} \binom{m}{a} \binom{m-a}{b} \binom{ta}{n-b}
(-1)^{m-a-b} (k-t)^b(z)^{n-m}.
$$
It follows that it must be the case that
$$
\sum_{m \geq n+1} \sum_{a=0}^m \sum_{b=0}^{m-a} \binom{m}{a}
\binom{m-a}{b} \binom{ta}{n-b} (-1)^{m-a-b} (k-t)^b(z)^{n-m} =0,
$$ so that
\begin{equation}\label{eq:4.4}
A_k = \sum_{n \geq 0} q^n \sum_{m =0}^n \sum_{a=0}^m
\sum_{b=0}^{m-a} (-1)^{m-a-b} \binom{m}{a} \binom{m-a}{b}
\binom{ta}{n-b} (k-t)^b(x_1-1)^{n-m}.
\end{equation}
Thus if we take the coefficient of $x_1^s$ on both sides of
(\ref{eq:4.4}) and we use the remark in the introduction that
$\overleftarrow{\des}_{[t]}(w)  = \overleftarrow{\ris}_{\{k+1-t,
\ldots, k\}}(w)$ for all $w \in [k]^n$, then we have the following
result.

\begin{theorem}\label{thm:4.1}
If $t \leq k$, then the number of words $w \in [k]^n$ such that
$\overleftarrow{\des}_{[t]}(w) =s$
{\rm(}$\overleftarrow{\ris}_{\{k+1-t, \ldots, k\}}(w) =s${\rm)} is
equal to
\begin{equation}\label{eq:4.5}
\sum_{m =0}^n \sum_{a=0}^m \sum_{b=0}^{m-a}
(-1)^{n-a-b-s} \binom{m}{a} \binom{m-a}{b} \binom{ta}{n-b} \binom{n-m}{s} (k-t)^b.
\end{equation}
\end{theorem}

If we want to find formulas for the number of words in $[k]^n$ with
$s$ descents that start with an element greater to $t$, then we need
to set $x_1 = y_1 = y_2 = z_1 = z_2 =1$ and $q_1 =q_2 =q$ in
(\ref{eq:4.1}). In that case, we will have $\lambda_j = 0$  and
$\nu_j = q$ for all $j$, $\mu_j = 1+q(x_2-1)$ for $j > t$, and
$\mu_j = 1$ for $j \leq  t$. It follows that
\begin{eqnarray*}
A_k &=& \frac{1}{1-\left( \sum_{j=0}^t q (1+q(x_2-1))^{k-t} +
\sum_{j=t+1}^k q\prod_{i=j+1}^k (1+q(x_2 -1)) \right)} \\
&=& \frac{1}{1-\left( qt(1+q(x_2-1))^{k-t} +
q\frac{(1+q(x_2-1))^{k-t} -1}{(1+q(x_2-1)) -1}\right)}\\
&=& \frac{1}{1-\frac{1}{(x_2-1)}\left( qt(x_2-1)(1+q(x_2-1))^{k-t} + (1+q(x_1-1))^{k-t} -1\right)}\\
&=& \frac{1}{1-\frac{1}{(x_2-1)}\left( (qt(x_2-1)+1)(1+q(x_2-1))^{k-t} -1\right)}\\
&=& \sum_{m \geq 0} \frac{1}{(x_2-1)^m} \left( (qt(x_2-1)+1)(1+q(x_2-1))^{k-t} -1\right)^m \\
&=& \sum_{m \geq 0} \frac{1}{(x_2-1)^m}\sum_{a=0}^m \binom{m}{a} (-1)^{m-a}
(qt(x_2-1)+1)^a (1+q(x_2-1))^{(k-t)a}\\
&=& \sum_{m \geq 0} \frac{1}{(x_2-1)^m}\sum_{a=0}^m \sum_{b=0}^{a}
\sum_{c=0}^{(k-t)a} (-1)^{m-a}\binom{m}{a} \binom{a}{b} \binom{(k-t)a}{c}
q^bt^b(x_2-1)^b q^c (x_2-1)^c.
\end{eqnarray*}
Again, if we want to take the coefficient of $q^n$, then we must have
$b+c =n$ or $c = n-b$. Thus
\begin{equation}\label{eq:4.6}
A_k = \sum_{n \geq 0} q^n \sum_{m \geq 0} \sum_{a=0}^m \sum_{b=0}^{a} (-1)^{m-a}\binom{m}{a} \binom{a}{b} \binom{(k-t)a}{n-b}
t^b (x_2-1)^n.
\end{equation}
Taking the coefficient of $q^n$ of both sides of
(\ref{eq:4.6}), we see that
\begin{equation}\label{eq:4.7}
\sum_{\pi \in [k]^n} x_2^{\overleftarrow{\des}_{\{t+1, \ldots k
\}}(\pi)} = \sum_{m \geq 0} \sum_{a=0}^m \sum_{b=0}^{a} (-1)^{m-a}
\binom{m}{a} \binom{a}{b} \binom{(k-t)a}{n-b} t^b (x_2-1)^{n-m}
\end{equation}
for all $n$. However, if we replace $x_2$ by $z+1$ in (\ref{eq:4.3}), we
see that the polynomial
$$\sum_{\pi \in [k]^n} (z+1)^{\overleftarrow{\des}_{\{t+1, \ldots, k \}}(\pi)}$$
has the Laurent expansion
$$
\sum_{m \geq 0} \sum_{a=0}^m \sum_{b=0}^{a} (-1)^{m-a}\binom{m}{a}
\binom{a}{b} \binom{(k-t)a}{n-b}t^b(z)^{n-m}.
$$
It follows that it must be the case that
$$
\sum_{m \geq n+1} \sum_{a=0}^m \sum_{b=0}^{a} (-1)^{m-a}\binom{m}{a}
\binom{a}{b} \binom{(k-t)a}{n-b}t^b(z)^{n-m} =0,
$$ so that
\begin{equation}\label{eq:4.8}
A_k = \sum_{n \geq 0} q^n \sum_{m =0}^n \sum_{a=0}^m \sum_{b=0}^{a}
(-1)^{m-a} \binom{m}{a} \binom{a}{b}
\binom{(k-t)a}{n-b}t^b(x_1-1)^{n-m}.
\end{equation}
Thus if we take the coefficient of $x_2^s$ on both sides of
(\ref{eq:4.8}) and we use the
remark in the introduction that $\overleftarrow{\des}_{\{t+1, \ldots, k\}}(w)  = \overleftarrow{\ris}_{[k-t]}(w)$
for all $w \in [k]^n$, then we have the following result.

\begin{theorem}\label{thm:4.2}
If $t \leq k$, then the number of words $w \in [k]^n$ such that
$\overleftarrow{\des}_{\{t+1, \ldots, k\}}(w)=s$
{\rm(}$\overleftarrow{\ris}_{[k-t]}(w) =s${\rm)} is equal to
\begin{equation}\label{eq:4.9}
\sum_{m =0}^n \sum_{a=0}^m \sum_{b=0}^{a}
(-1)^{n-a-s} \binom{m}{a} \binom{a}{b} \binom{(k-t)a}{n-b} \binom{n-m}{s} t^b.
\end{equation}
\end{theorem}

We end this section by showing that we can derive some non-trivial
binomial identities from Theorem \ref{thm:4.1} and \ref{thm:4.2}.
For example, in the special case of Theorem \ref{thm:4.2}
where $t=k-1$, we can count
the number of words $w \in [k]^n$ such that
$\overleftarrow{\des}_{\{k\}}(w)=s$ directly. We can classify the words in $[k]^n$ by how many $k$'s occur in the
word. That is, for those words $w \in [k]^n$ which have $n-r$ occurrences of
$k$, we can form a word such that $\overleftarrow{\des}_{\{k\}}(w)=s$ by first
picking a word $u \in [k-1]^{r}$. Next we insert a $k$ directly in front of
 $s$ different letters in $u$ in $\binom{r}{s}$ ways. Finally we can
place the remaining $n-r-s$ $k$'s either in a block with one of the
$k$'s that start a descent or at the end of $u$. The number ways to
place the remaining $k$'s is the number non-negative integer valued
solutions to $x_1 + \cdots + x_{s+1} = n-r-s$ or, equivalently, the
number of positive integer valued solutions to $y_1 + \cdots +
y_{s+1} = n-r+1$ which is clearly $\binom{n-r}{s}$. Note that to
have $s$ such descents, we must have $r \geq s$ and $n-r \geq s$ or,
equivalently, $s \leq r \leq n-s$. It follows that the number of
words $w \in [k]^n$ such that $\overleftarrow{\des}_{\{k\}}(w)=s$
equals
\begin{equation}\label{eq:4.10}
\sum_{r=s}^{n-s} (k-1)^{r} \binom{r}{s} \binom{n-r}{s}.
\end{equation}

Using (\ref{eq:4.9}) with $t =k-1$, we see that (\ref{eq:4.10}) equals
\begin{eqnarray}\label{eq:4.11}
&& \sum_{m =0}^n \sum_{a=0}^m \sum_{b=0}^{a}
(-1)^{n-a-s} \binom{m}{a} \binom{a}{b} \binom{a}{n-b} \binom{n-m}{s} (k-1)^b
= \nonumber \\
&& \sum_{b = 0}^n (k-1)^b \sum_{m=0}^n \sum_{a=b}^m
(-1)^{n-a-s} \binom{m}{a} \binom{a}{b} \binom{a}{n-b} \binom{n-m}{s}.
\end{eqnarray}
However, in (\ref{eq:4.11}), we must have $n-m \geq s$ or, equivalently,
$n-s \geq m$. Since $m \geq a \geq b$, we must have $b \geq s$ since
otherwise the binomial coefficient $\binom{a}{n-b}$ will equal 0.
Thus (\ref{eq:4.10}) equals
\begin{equation}\label{eq:4.12}
\sum_{b = s}^{n-s} (k-1)^b \sum_{m=b}^{n-s} \sum_{a=b}^m
(-1)^{n-a-s} \binom{m}{a} \binom{a}{b} \binom{a}{n-b} \binom{n-m}{s}.
\end{equation}
Since (\ref{eq:4.10}) and (\ref{eq:4.12}) hold for all $k$, it
follows that
\begin{equation}\label{eq:4.13}
\sum_{r=s}^{n-s} x^{r} \binom{r}{s} \binom{n-r}{s} =
\sum_{b=s}^{n-s} x^b \sum_{m=0}^n \sum_{a=b}^m
(-1)^{n-a-s} \binom{m}{a} \binom{a}{b} \binom{a}{n-b} \binom{n-m}{s}.
\end{equation}
Thus we have proved that the following identity holds.
\begin{equation}\label{eq:4.14}
\binom{r}{s} \binom{n-r}{s} =
\sum_{m=r}^{n-s} \sum_{a=r}^m
(-1)^{n-a-s} \binom{m}{a} \binom{a}{r} \binom{a}{n-r} \binom{n-m}{s}.
\end{equation}

Next consider the special case of Theorem \ref{thm:4.1} where
$t =2$.  Suppose we want to count the number of $\pi \in [k]^n$ where
$\overleftarrow{\des}_{[2]}(w)=s$. We can classify such words according to
number $a$ of $1$'s and the number $b$ of 2's that appear in the word.
Clearly since the only descents that we can count are cases where there
is a 2 followed by a 1, we must have $a,b \geq s$. We claim that
we can count such words as follows. First we pick
word $w$ of length $n -a -b$ made up of letters from $\{3, \ldots, k\}$ in
$(k-2)^{n-a-b}$ ways. Then to create the $s$ $2~1$ descents, we imagine
inserting letters of the form $\overline{21}$ into $w$ to get a word
$u$ of length $n-a-b+s$ over the alphabet $\{\overline{21},3, \ldots,k\}$. This can be done in  $\binom{n-a-b +s}{s}$ ways. For each
such $u$, we first insert the remaining $a-s$ 1's to get word $v$ of
length $n-a-b+s +a -s = n-b$ over the alphabet
$\{1,\overline{21},3, \ldots,k\}$. Since we can
insert the 1's in front of any of the letters of $u$ over the
the alphabet $\{\overline{21},3, \ldots,k\}$ or at the end, the number of
ways to insert the remaining 1's is the number of nonnegative integer solutions
to $x_1+ \cdots + x_{n-a-b+s+1} = a-s$ which is $\binom{n-b}{a-s}$. Finally, we
have to insert the remaining $b-s$ 2's. In this case, since we can insert the
2's into $v$ in front of any letter which is not a 1 or at the end, the number of ways to insert the remaining 2's is the number of nonnegative integer solutions to $x_1+ \cdots +x_{n-a-b+s+1} = b-s$ which is $\binom{n-a}{b-s}$. Thus it
follows that the number of words $\pi \in [k]^n$ such that
$\overleftarrow{\des}_{[2]}(w)=s$ is
\begin{equation}\label{eq:4.15}
\sum_{\overset{a,b\geq s}{a+b \leq n}} (k-2)^{n-a-b} \binom{n-a
-b+s}{s} \binom{n-b}{a-s}\binom{n-a}{b-s}.
\end{equation}
On the other hand, from Theorem \ref{thm:4.1}, we see that the
number of words $\pi \in [k]^n$ such that
$\overleftarrow{\des}_{[2]}(w)=s$ is
\begin{equation}\label{eq:4.16}
\sum_{m =0}^n \sum_{a=0}^m \sum_{b=0}^{m-a}
(-1)^{n-a-b-s} \binom{m}{a} \binom{m-a}{b} \binom{2a}{n-b} \binom{n-m}{s} (k-2)^b.
\end{equation}
Since (\ref{eq:4.15}) and (\ref{eq:4.16}) hold for all $k$, it must
be the case that
\begin{eqnarray}\label{eq:4.17}
&&\sum_{\overset{a,b\geq s}{a+b \leq n}} x^{n-a-b} \binom{n-a -b+s}{s}
\binom{n-b}{a-s}\binom{n-a}{b-s}  \nonumber \\
&&\quad=\sum_{m =0}^n \sum_{a=0}^m \sum_{b=0}^{m-a} (-1)^{n-a-b-s}
\binom{m}{a} \binom{m-a}{b} \binom{2a}{n-b} \binom{n-m}{s} x^b.
\end{eqnarray}
Taking the coefficient of $x^r$ on both sides yields the following identity.
\begin{equation}\label{eq:4.18}
\sum_{a=s}^{n-s} \binom{r+s}{s} \binom{a+r}{a-s} \binom{n-a}{n-a -r -s} =
\sum_{m=0}^n\sum_{a=0}^{m-r} (-1)^{n-a-r-s}
\binom{m}{a} \binom{m-a}{r} \binom{2a}{n-r} \binom{n-m}{s}.
\end{equation}

\section{Classifying descents and rises by their equivalence
classes $\mod s$ for $s \geq 2$.}\label{sec5}

In this section we study the set partition
$\mathbb{N}=\mathbb{N}_1\cup\mathbb{N}_2
\cup \cdots \cup \mathbb{N}_s$ where $s > 2$ and $\mathbb{N}_i=\{j\ |\ j= i\mod
s\}$ for $i =1, \ldots, s$. In this case, we shall denote
$N_i = s\mathbb{N}+i$ for $i =1, \ldots, s-1$ and $N_s = s\mathbb{N}$.

Recall that we can rewrite (\ref{eq:Main}) as
\begin{equation}\label{eq:Main22}
A_k= \frac{1+ \sum_{j=1}^k \lambda_j \prod_{i=j+1}^k \mu_i}{1- \sum_{j=1}^k \nu_j \prod_{i=j+1}^k \mu_i}
\end{equation}
where $\lambda_j = \frac{q_j(1-y_j)}{1-q_j(z_j-y_j)}$,
$\mu_i = \frac{q_i(z_i-x_i)}{1-q_i(z_i-y_i)}$, and
$\nu_j = \frac{q_jy_j}{1-q_j(z_j-y_j)}$.

We let $A^{(s)}_k$ denote $A_k$ under the substitution that
$\lambda_{si+j} = \lambda_j$, $\mu_{si+j} = \mu_j$, and $\nu_{si+j} = \nu_j$  for all $i$ and  $j=1,\ldots, s$.
Then it is easy to see that for $k \geq 1$,
\begin{eqnarray}\label{eq:Ask}
&\ \ \ A^{(s)}_{sk} &= \frac{1+\left(\sum_{j=1}^s \lambda_j
\prod_{i=j+1} \mu_i\right)
\left(\sum_{r=0}^{k-1} (\mu_1 \mu_2 \cdots \mu_s)^r\right)}{1-\left(\sum_{j=1}^s \nu_j \prod_{i=j+1} \mu_i\right)\left(\sum_{r=0}^{k-1} (\mu_1 \mu_2 \cdots \mu_s)^r\right)} \nonumber \\
&&=\frac{1+\left(\sum_{j=1}^s \lambda_j \prod_{i=j+1} \mu_i\right)
\left(\frac{(\mu_1 \mu_2 \cdots \mu_s)^k-1}{(\mu_1 \mu_2 \cdots
\mu_s)-1} \right)}{1-\left(\sum_{j=1}^s \nu_j \prod_{i=j+1}
\mu_i\right) \left(\frac{(\mu_1 \mu_2 \cdots \mu_s)^k-1}{(\mu_1
\mu_2 \cdots \mu_s)-1} \right)}.\nonumber
\end{eqnarray}

More generally, we can express $A^{(s)}_k$ in
the form
$$A^{(s)}_k =  \frac{\Theta(\lambda_1, \ldots, \lambda_k,\mu_1, \ldots, \mu_k)}{\Theta(-\nu_1, \ldots, -\nu_k,\mu_1, \ldots, \mu_k)}$$
where
\begin{equation}\label{theta}
\Theta(\lambda_1, \ldots, \lambda_k,\mu_1, \ldots, \mu_k) = 1+
\left(\sum_{j=1}^k \lambda_j \prod_{i=j+1}^k \mu_i\right)
\left(\sum_{r=0}^{k-1} (\mu_1 \mu_2 \cdots \mu_s)^r\right).
\end{equation}
Then for $1 \leq t < s$, we have that
\begin{eqnarray}
&& \ \Theta(\lambda_1, \ldots, \lambda_k,\mu_1, \ldots, \mu_k) \\
&&=1+ (\sum_{j=1}^t \lambda_j\prod_{i=j+1}^t \mu_i) (\sum_{r=0}^k
(\mu_1 \mu_2 \cdots \mu_s)^r) + (\mu_1 \mu_2 \cdots \mu_t)
(\sum_{j=1}^{t+1} \lambda_j\prod_{i=j+1}^s \mu_i)(\sum_{r=0}^{k-1} (\mu_1 \mu_2 \cdots \mu_s)^r)  \nonumber \\
&&=1+ (\sum_{j=1}^t \lambda_j\prod_{i=j+1}^t
\mu_i)\left(\frac{(\mu_1 \mu_2 \cdots \mu_s)^{k+1}-1}{(\mu_1 \mu_2
\cdots \mu_s)-1}\right) + (\mu_1 \mu_2 \cdots \mu_t)
(\sum_{j=1}^{t+1} \lambda_j\prod_{i=j+1}^s \mu_i)\left(\frac{(\mu_1
\mu_2 \cdots \mu_s)^{k}-1}{(\mu_1 \mu_2 \cdots \mu_s)-1}\right).
\nonumber
\end{eqnarray}
Hence, for $1 \leq t \leq s$,
\begin{eqnarray}\label{eq:Ask+t}
&& \ \ \ A^{(s)}_{sk+t} = \\
&&\frac{1+ (\sum_{j=1}^t \lambda_j\prod_{i=j+1}^t \mu_i)\left(\frac{(\mu_1 \mu_2 \cdots \mu_s)^{k+1}-1}{(\mu_1 \mu_2 \cdots \mu_s)-1}\right) + (\mu_1 \mu_2 \cdots \mu_t)
(\sum_{j=1}^{t+1} \lambda_j\prod_{i=j+1}^s \mu_i)\left(\frac{(\mu_1 \mu_2 \cdots \mu_s)^{k}-1}{(\mu_1 \mu_2 \cdots \mu_s)-1}\right)}{1- (\sum_{j=1}^t \nu_j\prod_{i=j+1}^t \mu_i)\left(\frac{(\mu_1 \mu_2 \cdots \mu_s)^{k+1}-1}{(\mu_1 \mu_2 \cdots \mu_s)-1}\right) - (\mu_1 \mu_2 \cdots \mu_t)
(\sum_{j=1}^{t+1} \nu_j\prod_{i=j+1}^s \mu_i)\left(\frac{(\mu_1 \mu_2 \cdots \mu_s)^{k}-1}{(\mu_1 \mu_2 \cdots \mu_s)-1}\right)}. \nonumber
\end{eqnarray}

\subsection{The case where $k$ is equal to 0 mod $s$.}

First we shall consider formulas for the number of words in $[sk]^n$ with
$p$ descents whose first element is equivalent to $r \mod s$ where $1 \leq r \leq s$. Note that if we consider the complement map
$comp_{sk}:[sk]^n \rightarrow
[sk]^n$ given by $comp(\pi_1 \cdots \pi_n) = (sk+1-\pi_1)\cdots (sk+1-\pi_n)$,
then it is easy to see that
$\overleftarrow{\des}_{s\mathbb{N}+r}(\pi) =
\overleftarrow{\ris}_{s\mathbb{N}+s+1-r}(comp_{sk}(\pi))$ for
$r =1, \ldots, s$. Thus the problem of counting the number of
words in $[sk]^n$  with $p$ descents whose first element is equivalent to $r \mod s$ is the same as counting the number of
words in $[sk]^n$  with $p$ rises whose first element is equivalent to
$s+1 -r \mod s$

Now consider the case where $z_i =y_i =1$ and $q_i=q$ for
$i =1, \ldots, s$ and $x_i =1$ for $i \neq r$. In this case,
$$A^{(s)}_{sk} = \sum_{n\geq 0} q^n \sum_{\pi \in [sk]^n} x_r^{\overleftarrow{\des}_{s\mathbb{N}+r}(\pi)}.$$
Substituting into our formulas for $A^{(s)}_{sk}$, we see that in
this case $\lambda_i =0$ and $\nu_i =q$ for $i=1, \ldots, s$ and
$\mu_i=1$ for $i\neq r$ and $\mu_r = 1+q(x_r-1)$. Thus
under this substitution, (\ref{eq:Ask}) becomes
\begin{eqnarray}\label{skr:1}
&&\ \ \ A^{(s)}_{sk} \\
&&= \frac{1}{1-((r-1)q\mu_r+(s-r+1)q)\frac{\mu_r^k -1}{q(x_r-1)}}  \nonumber \\
&&= \frac{1}{1-\frac{1}{(x_r-1)} (s+(s-1)q(x_r-1))(\mu_r^k -1)} \nonumber\\
&&= \sum_{j=0}^\infty \frac{1}{(x_r-1)^j} (s+(r-1)q(x_r-1))^j
\mu_r^k -1)^j  \nonumber \\
&&= \sum_{j=0}^\infty \frac{1}{(x_r-1)^j} \sum_{i_1,i_2 = 0}^j
\binom{j}{i_1} s^{j-i_1}(r-1)^{i_1}q^{i_1}(x_r-1)^{i_1}
\binom{j}{i_2}(-1)^{j-i_2}\mu_r^{ki_2}  \nonumber \\
&&= \sum_{j=0}^\infty \frac{1}{(x_r-1)^j}\sum_{i_1,i_2 = 0}^j
\binom{j}{i_1}s^{j-i_1}(r-1)^{i_1}q^{i_1}(x_r-1)^{i_1}
\binom{j}{i_2}(-1)^{j-i_2} \left(\sum_{t=0}^{ki_2} \binom{ki_2}{t}
q^t (x_r-1)^t\right). \nonumber
\end{eqnarray}
Taking the coefficient of $q^n$ in (\ref{skr:1}), we see that
$n =t+i_1$ so that
\begin{equation}\label{skr:2}
A^{(s)}_{sk} =
\sum_{n \geq 0} q^n
\sum_{j=0}^\infty \sum_{i_1,i_2 = 0}^j (-1)^{j-i_2}
s^{j-i_1}(r-1)^{i_1}\binom{j}{i_1}
\binom{j}{i_2} \binom{ki_2}{n-i_1}  (x_r-1)^{n-j}.
\end{equation}
Thus we must have
\begin{equation}\label{skr:3}
\sum_{\pi \in [sk]^n} x_r^{\overleftarrow{\des}_{s\mathbb{N}+r}(\pi)} =
 \sum_{j=0}^\infty \sum_{i_1,i_2 = 0}^j (-1)^{j-i_2}
s^{j-i_1}(r-1)^{i_1}\binom{j}{i_1}
\binom{j}{i_2} \binom{ki_2}{n-i_1}  (x_r-1)^{n-j}
\end{equation}
for all $n$. However, if we replace $x_r$ by $z+1$ in (\ref{skr:3}), we
see that the polynomial
$$\sum_{\pi \in [sk]^n} (z+1)^{\overleftarrow{\des}_{s\mathbb{N}+r}(\pi)}$$
has the Laurent expansion
$$\sum_{j=0}^\infty \sum_{i_1,i_2 = 0}^j (-1)^{j-i_2}
s^{j-i_1}(r-1)^{i_1}\binom{j}{i_1}
\binom{j}{i_2} \binom{ki_2}{n-i_1}  (z)^{n-j}.$$
It follows that it must be the case that
$$\sum_{j=n+1}^\infty \sum_{i_1,i_2 = 0}^j (-1)^{j-i_2}
s^{j-i_1}(s-1)^{i_1}\binom{j}{i_1}
\binom{j}{i_2} \binom{ki_2}{n-i_1}  (z)^{n-j} = 0$$ so that
\begin{equation}\label{skr:4}
A^{(s)}_{sk} =
\sum_{n \geq 0} q^n
\sum_{j=0}^n \sum_{i_1,i_2 = 0}^j (-1)^{j-i_2}
s^{j-i_1}(s-1)^{i_1}\binom{j}{i_1}
\binom{j}{i_2} \binom{ki_2}{n-i_1}  (x_r-1)^{n-j}.
\end{equation}

Thus we have the following theorem by taking the coefficient of
$x_r^p$ on both sides of (\ref{skr:4}).
\begin{theorem}\label{thm:skr}
The number of words $\pi \in [sk]^n$ with
$\overleftarrow{\des}_{s\mathbb{N}+r}(\pi) =p$
($\overleftarrow{\ris}_{s\mathbb{N}+s+1 -r}(\pi) =p$ )
is
\begin{equation}\label{eq:thm:skr}
\sum_{j=0}^n \sum_{i_1,i_2 = 0}^j (-1)^{n+p+i_2}
s^{j-i_1}(s-1)^{i_1} \binom{j}{i_1}
\binom{j}{i_2} \binom{ki_2}{n-i_1}  \binom{n-j}{p}.
\end{equation}
\end{theorem}

In the case $s=2$, our formulas simplify somewhat. For example,
putting $s =2$ and $r=2$ in Theorem \ref{thm:skr}, we obtain the following.
\begin{corollary}\label{even2k}
The number of $n$-letter words $\pi$ on $[2k]$ having
$\overleftarrow{\des}_E(\pi)=p$ {\rm(}resp.
$\overleftarrow{\ris}_O(\pi)=p${\rm)} is given by
$$\sum_{j=0}^n\sum_{i_1,i_2=0}^j(-1)^{n+p+i_2}2^{j-i_1}\binom{j}{i_1}\binom{j}{i_2}\binom{ki_2}{n-i_1}\binom{n-j}{p}.$$
\end{corollary}

Similarly, putting $s =2$ and $r=1$ in Theorem \ref{thm:skr}, we obtain the following.
\begin{corollary}
The number of $n$-letter words $\pi$ on $[2k]$ having
$\overleftarrow{\des}_O(\pi)=p$ {\rm(}resp.
$\overleftarrow{\ris}_E(\pi)=p${\rm)} is given by
$$\sum_{j=0}^n\sum_{i=0}^j(-1)^{n+p+i}2^j\binom{j}{i}\binom{i}{n}\binom{n-j}{p}.$$
\end{corollary}

\subsection{The cases where $k$ is equal to $t$ mod $s$ for $t=1,\ldots,s-1$.}

Fix $t$ where $1 \leq t \leq s-1$. First we shall consider formulas for the number of words in $[sk+t]^n$ with
$p$ descents whose first element is equivalent to $r \mod s$ where $1 \leq r \leq s$. We shall see that we have to divide this problem into two cases depending on whether $r \leq t$ or $r > t$.
Note that if we consider the complement map
$comp_{sk+t}:[sk+t]^n \rightarrow
[sk+t]^n$ given by $comp(\pi_1 \cdots \pi_n) = (sk+t+1-\pi_1)\cdots (sk+t+1-\pi_n)$,
then it is easy to see that
$\overleftarrow{\des}_{s\mathbb{N}+r}(\pi) =
\overleftarrow{\ris}_{s\mathbb{N}+t+1-r}(comp_{sk}(\pi))$ for
$r =1, \ldots, t$ and $\overleftarrow{\des}_{s\mathbb{N}+r}(\pi) =
\overleftarrow{\ris}_{s\mathbb{N}+s+r-t-1}(comp_{sk}(\pi))$ for
$r =t+1, \ldots, s$.

First consider the case where $y_i=z_i =1$ for $i =1, \ldots, s$ and
$x_i=1$ for $i \neq r$ where $r >t$. In this case,
$$A^{(s)}_{sk+t} = \sum_{n\geq 0} q^n \sum_{\pi \in [sk+t]^n} x_r^{\overleftarrow{\des}_{s\mathbb{N}+r}(\pi)}.$$
Substituting into our formulas for $A^{(s)}_{sk+t}$, we see that in
this case $\lambda_i =0$ and $\nu_i =q$ for $i=1, \ldots, s$ and
$\mu_i=1$ for $i\neq r$ and $\mu_r = 1+q(x_r-1)$. Thus
under this substitution, (\ref{eq:Ask}) becomes
\begin{eqnarray}\label{sktr:1}
&&\ \ \ A^{(s)}_{sk+t} \\
&&= \frac{1}{1-qt\frac{\mu_r^{k+1} -1}{q(x_r-1)}-
((r-1-t)q\mu_r +(s-r+1)q) \frac{\mu_r^{k} -1}{q(x_r-1)}} \nonumber \\
&&= \frac{1}{1-\frac{1}{(x_r-1)}[t\mu_r)^{k+1} -1)+
(s-t+ (r-1-t)q(x_r-1))(\mu_r^{k} -1)]} \nonumber \\
&&= \frac{1}{1-\frac{1}{(x_r-1)}[\mu_r^{k}[s + (r-1)q(x_r-1)]-
[s+(r-1-t)q(x_r-1)]]}\nonumber \\
&&=\sum_{m=0}^\infty \frac{1}{(x_r-1)^m} [\mu_r^{k}[s +
(r-1)q(x_r-1)]-
[s+(r-1-t)q(x_r-1)]]^m \nonumber \\
&&=\sum_{m=0}^\infty  \sum_{j=0}^m \frac{(-1)^{m-j}}{(x_r-1)^m}
\binom{m}{j} (s+(r-1-t)q(x_r-1))^{m-j} (s + (r-1)q(x_r-1))^j
\mu_r^{kj}. \nonumber
\end{eqnarray}
Using the expansions
\begin{eqnarray*}
(s+(r-1-t)q(x_r-1))^{m-j} &=& \sum_{i_1 =0}^{m-j} \binom{m-j}{i_1} s^{m-j-i_1}
(r-1-t)^{i_1} q^{i_1} (x_r-1)^{i_1}, \\
 (s + (r-1)q(x_r-1))^j &=& \sum_{i_2 =0}^{j} \binom{j}{i_2} s^{j-i_2}
(r-1)^{i_2} q^{i_2} (x_r-1)^{i_2}, \ \mbox{and}  \\
\mu_r^{kj} &=& \sum_{i_3 =0}^{kj} \binom{kj}{i_3} q^{i_3}
(x_r-1)^{i_3},
\end{eqnarray*}
and setting $i_1 + i_2 +i_3 =n$, we see that (\ref{sktr:1}) becomes
\begin{eqnarray}\label{sktr:2}
&&\ \ \ A^{(s)}_{sk+t} =\\
&&\sum_{n \geq 0} q^n \sum_{m=0}^\infty \sum_{i_1 =0}^{m-j}
\sum_{i_2 =0}^{j}(-1)^{m-j} s^{m-i_1-i_2} (r-1-t)^{i_1} (r-1)^{i_2}
\binom{m}{j} \binom{m-j}{i_1} \binom{j}{i_2} \binom{kj}{n-i_1-i_2}
(x_r-1)^{n-m}. \nonumber
\end{eqnarray}
Thus we must have
\begin{eqnarray}\label{sktr:3}
&&\sum_{\pi \in [sk_t]^n} x_r^{\overleftarrow{\des}_{s\mathbb{N}+r}(\pi)}  \\
&&=\sum_{m=0}^\infty \sum_{i_1 =0}^{m-j} \sum_{i_2 =0}^{j}(-1)^{m-j}
s^{m-i_1-i_2} (r-1-t)^{i_1} (r-1)^{i_2} \binom{m}{j}
\binom{m-j}{i_1} \binom{j}{i_2} \binom{kj}{n-i_1-i_2} (x_r-1)^{n-m}.
\nonumber
\end{eqnarray}
for all $n$. However, if we replace $x_r$ by $z+1$ in (\ref{sktr:3}), we
that the polynomial
$$\sum_{\pi \in [sk+t]^n} (z+1)^{\overleftarrow{\des}_{s\mathbb{N}+r}(\pi)}$$
has the Laurent expansion
$$\sum_{m=0}^\infty \sum_{i_1 =0}^{m-j}
\sum_{i_2 =0}^{j}(-1)^{m-j} s^{m-i_1-i_2} (r-1-t)^{i_1} (r-1)^{i_2}
\binom{m}{j} \binom{m-j}{i_1} \binom{j}{i_2} \binom{kj}{n-i_1-i_2}
(z)^{n-m}. $$
It follows that it must be the case that
$$\sum_{m=n+1}^\infty \sum_{i_1 =0}^{m-j}
\sum_{i_2 =0}^{j}(-1)^m-j s^{m-i_1-i_2} (r-1-t)^{i_1} (r-1)^{i_2}
\binom{m}{j} \binom{m-j}{i_1} \binom{j}{i_2} \binom{kj}{n-i_1-i_2}
(x_r-1)^{n-m} = 0$$ so that
\begin{eqnarray}\label{sktr:4}
&& \ \ \ \ \ \ A^{(s)}_{sk+t}  \\
&&=\sum_{m=0}^n \sum_{i_1 =0}^{m-j} \sum_{i_2 =0}^{j}(-1)^{m-j}
s^{m-i_1-i_2} (r-1-t)^{i_1} (r-1)^{i_2} \binom{m}{j}
\binom{m-j}{i_1} \binom{j}{i_2} \binom{kj}{n-i_1-i_2} (x_r-1)^{n-m}.
\nonumber
\end{eqnarray}

Thus we have the following theorem by taking the coefficient of
$x_r^p$ on both sides of (\ref{sktr:4}).
\begin{theorem}\label{thm:sktr} If $t =1, \ldots, s-1$ and
$t<r \leq s$, then
the number of words $\pi \in [sk+t]^n$ with
$\overleftarrow{\des}_{s\mathbb{N}+r}(\pi) =p$
($\overleftarrow{\ris}_{s\mathbb{N}+s+r-t-1}(\pi) =p$ )
is
\begin{equation}\label{eq:thm:sktr}
\sum_{m=0}^n \sum_{i_1 =0}^{m-j}
\sum_{i_2 =0}^{j}(-1)^{n+p+j} s^{m-i_1-i_2} (r-1-t)^{i_1} (r-1)^{i_2}
\binom{m}{j} \binom{m-j}{i_1} \binom{j}{i_2} \binom{kj}{n-i_1-i_2}
\binom{n-m}{p}
\end{equation}
\end{theorem}

In the case $s=2$, our formulas simplify somewhat. For example,
putting $s =2$, $r=2$ and $t=1$ in Theorem \ref{thm:sktr}, we obtain the following.

\begin{corollary}
The number of $n$-letter words $\pi$ over $[2k+1]$ having
$\overleftarrow{\des}_E(\pi)=p$ {\rm(}resp.
$\overleftarrow{\ris}_E(\pi)=p${\rm)} is given by
$$\sum_{m=0}^n\sum_{j=0}^m\sum_{i=0}^j(-1)^{n+p+j}2^{m-i}\binom{m}{j}\binom{j}{i}\binom{kj}{n-i}\binom{n-m}{p}.$$
\end{corollary}

Next consider the case where $y_i=z_i =1$ for $i =1, \ldots, s$ and
$x_i=1$ for $i \neq r$ where $r \leq t$. In this case,
$$A^{(s)}_{sk+t} = \sum_{n\geq 0} q^n \sum_{\pi \in [sk+t]^n} x_r^{\overleftarrow{\des}_{s\mathbb{N}+r}(\pi)}.$$
Substituting into our formulas for $A^{(s)}_{sk+t}$, we see that in
this case $\lambda_i =0$ and $\nu_i =q$ for $i=1, \ldots, s$ and
$\mu_i=1$ for $i\neq r$ and $\mu_r = 1+q(x_r-1)$. Thus
under this substitution, (\ref{eq:Ask}) becomes
\begin{eqnarray}\label{1sktr:1}
&&\ \ \ \ \ \ \ \ \ A^{(s)}_{sk+t} \\
&&= \frac{1}{1-((r-1)q\mu_r) + q(t-r+1))\frac{\mu_r^{k+1} -1}{q(x_r-1)}- (s-t)q \mu_r\frac{\mu_r^{k} -1}{q(x_r-1)}} = \nonumber \\
&& \frac{1}{1-\frac{1}{(x_r-1)}[(t +(r-1)q(x_r-1))(\mu_r^{k+1} -1)
+ (s-t)\mu_r(\mu_r^{k} -1)]} \nonumber \\
&&= \frac{1}{1-\frac{1}{(x_r-1)}[\mu_r^{k+1}[s+(r-1)q(x_r-1)]-
[s+(s-t+r-1)q(x_r-1)]]} \nonumber \\
&&=\sum_{m=0}^\infty \frac{1}{(x_r-1)^m} [\mu_r^{k+1}[s +
(r-1)q(x_r-1)]-
[s+(s-t+r-1)q(x_r-1)]]^m  \nonumber \\
&&=\sum_{m=0}^\infty \sum_{j=0}^m  \frac{(-1)^{m-j}}{(x_r-1)^m}
\binom{m}{j} (s+(s-t+r-1)q(x_r-1))^{m-j} (s + (r-1)q(x_r-1))^j
\mu_r^{kj+j}. \nonumber
\end{eqnarray}
Using the expansions
\begin{eqnarray*}
(s+(s-t+r-1)q(x_r-1))^{m-j} &=& \sum_{i_1 =0}^{m-j} \binom{m-j}{i_1} s^{m-j-i_1}(s-t+r-1)^{i_1} q^{i_1} (x_r-1)^{i_1}, \\
 (s + (r-1)\mu_r^j &=& \sum_{i_2 =0}^{j} \binom{j}{i_2} s^{j-i_2}
(r-1)^{i_2} q^{i_2} (x_r-1)^{i_2}, \ \mbox{and}  \\
\mu_r^{kj+j} &=& \sum_{i_3 =0}^{kj+j} \binom{kj+j}{i_3} q^{i_3}
(x_r-1)^{i_3},
\end{eqnarray*}
and setting $i_1 + i_2 +i_3 =n$, we see that (\ref{1sktr:1}) becomes
\begin{eqnarray}\label{1sktr:2}
&&\ \ \ \ \ \ \ \ A^{(s)}_{sk+t} =\\
&&\sum_{n \geq 0} q^n \sum_{m=0}^\infty \sum_{i_1 =0}^{m-j}
\sum_{i_2 =0}^{j}(-1)^{m-j} s^{m-i_1-i_2} (s-t+r-1)^{i_1} (r-1)^{i_2}
\binom{m}{j} \binom{m-j}{i_1} \binom{j}{i_2} \binom{kj+j}{n-i_1-i_2}
(x_r-1)^{n-m}. \nonumber
\end{eqnarray}
Thus we must have
\begin{eqnarray}\label{1sktr:3}
&&\sum_{\pi \in [sk_t]^n} x_r^{\overleftarrow{\des}_{s\mathbb{N}+r}(\pi)} = \\
&&\sum_{m=0}^\infty \sum_{i_1 =0}^{m-j}
\sum_{i_2 =0}^{j}(-1)^{m-j} s^{m-i_1-i_2} (s-t+r-1)^{i_1} (r-1)^{i_2}
\binom{m}{j} \binom{m-j}{i_1} \binom{j}{i_2} \binom{kj+j}{n-i_1-i_2}
(x_r-1)^{n-m}. \nonumber
\end{eqnarray}
for all $n$. However, if we replace $x_r$ by $z+1$ in (\ref{1sktr:3}), we
that the polynomial
$$\sum_{\pi \in [sk+t]^n} (z+1)^{\overleftarrow{\des}_{s\mathbb{N}+r}(\pi)}$$
has the Laurent expansion
$$\sum_{m=0}^\infty \sum_{i_1 =0}^{m-j}
\sum_{i_2 =0}^{j}(-1)^m-j s^{m-i_1-i_2} (s-t+r-1)^{i_1} (r-1)^{i_2}
\binom{m}{j} \binom{m-j}{i_1} \binom{j}{i_2} \binom{kj+j}{n-i_1-i_2}
(z)^{n-m}. $$
It follows that it must be the case that
$$\sum_{m=n+1}^\infty \sum_{i_1 =0}^{m-j}
\sum_{i_2 =0}^{j}(-1)^{m-j} s^{m-i_1-i_2} (s-t+r-1)^{i_1} (r-1)^{i_2}
\binom{m}{j} \binom{m-j}{i_1} \binom{j}{i_2} \binom{kj+j}{n-i_1-i_2}
(x_r-1)^{n-m} = 0$$ so that
\begin{eqnarray}\label{1sktr:4}
&& \ \ \ \ \ \ A^{(s)}_{sk+t} = \\
&&\sum_{m=0}^\infty \sum_{i_1 =0}^{m-j}
\sum_{i_2 =0}^{j}(-1)^{m-j} s^{m-i_1-i_2} (s-t+r-1)^{i_1} (r-1)^{i_2}
\binom{m}{j} \binom{m-j}{i_1} \binom{j}{i_2} \binom{kj+j}{n-i_1-i_2}
(x_r-1)^{n-m}. \nonumber
\end{eqnarray}

Thus we have the following theorem by taking the coefficient of
$x_r^p$ on both sides of (\ref{sktr:4}).
\begin{theorem}\label{thm:1sktr} If $k \geq 0$, $s \geq 2$,
$t =1, \ldots, s-1$,  and
$t<r \leq s$, then
the number of words $\pi \in [sk+t]^n$ with
$\overleftarrow{\des}_{s\mathbb{N}+r}(\pi) =p$
($\overleftarrow{\ris}_{s\mathbb{N}+s+r-t-1}(\pi) =p$ )
is
\begin{equation}\label{eq:thm:1sktr}
\sum_{m=0}^\infty \sum_{i_1 =0}^{m-j}
\sum_{i_2 =0}^{j}(-1)^{n+p+j} s^{m-i_1-i_2} (s-t+r-1)^{i_1} (r-1)^{i_2}
\binom{m}{j} \binom{m-j}{i_1} \binom{j}{i_2} \binom{kj+j}{n-i_1-i_2}
\binom{n-m}{p}
\end{equation}
\end{theorem}

In the case $s=2$, our formulas simplify somewhat. For example,
putting $s =2$, $r=1$ and $t=1$ in Theorem \ref{thm:sktr}, we obtain the following.

\begin{corollary}
The number of $n$-letter words $\pi$ over $[2k+1]$ having
$\overleftarrow{\des}_O(\pi)=p$ {\rm(}resp.
$\overleftarrow{\ris}_O(\pi)=p${\rm)} is given by
$$\sum_{m=0}^n\sum_{j=0}^m\sum_{i=0}^{kj+j}(-1)^{n+p+j}2^{m-n+i}\binom{m}{j}\binom{m-j}{n-i}\binom{jk+j}{i}\binom{n-m}{p}.$$
\end{corollary}

\section{Concluding remarks}\label{sec6}

A particular case of the results obtained by Burstein and Mansour
in~\cite{BM} is the distribution of descents (resp. levels, rises),
which can be viewed as occurrences of so called {\em generalized
patterns} 21 (resp. 11, 12) in words. To get these distributions
from our results, we proceed as follows (we explain only the case of
descents; rises and levels can be considered similarly). Set
$x_1=x_2=x$, $y_1=y_2=z_1=z_2=1$, and $q_1=q_2=q$ in $A^{(2)}_{2k}$
and $A^{(2)}_{2k+1}$ to get the distribution in \cite[Theorem
2.2]{BM} for $\ell=2$ (the case of descents/rises). Thus, our
results refine and generalize the known
distributions of descents, levels, and rises in words. \\

It is interesting to compare our formulas with formulas of Hall and
Remmel \cite{HR}.  For example, suppose that $X=E$ and $Y = \mathbb{N}$ and
$\rho =(\rho_1, \ldots, \rho_{2k})$ is a composition of $n$. Then
Theorem \ref{combword1} tells that the number of words $\pi$ of $[2k]^n$ such
that $\overleftarrow{\des}_E(\pi) =p$ is
\begin{equation}\label{HRE}
\binom{a}{\rho_{2},\rho_{4}, \ldots, \rho_{2k}}
\sum_{r=0}^p (-1)^{p-r} \binom{a+r}{r}
\binom{n+1}{p-r}
\prod_{i=1}^k \binom{\rho_{2i} + r +
(\rho_{2i+1}+ \rho_{2i+3}+ \cdots +\rho_{2k-1})}{\rho_{2i}},
\end{equation}
where $a =\rho_2+\rho_4 + \cdots +\rho_{2k}$. This shows that once
we are given the distribution of the letters for words in $[2k]^n$,
we can find an expression for the number of words $\pi$ such that
$\overleftarrow{\des}_E(\pi) =p$ with a single alternating sum of
products of binomial coefficients. This contrasts with Corollary
\ref{even2k} where we require a triple alternating sum of products
of binomial coefficients to get an expression for the number of
words of $[2k]^n$ such that $\overleftarrow{\des}_E(\pi) =p$. Of
course, we can get a similar expressions for the number of words of
$[2k]^n$ such that $\overleftarrow{\des}_E(\pi) =p$ by summing the
formula in (\ref{HRE}) over all $\binom{n+k-1}{k-1}$ compositions of
$n$ into $k$ parts but that has the disadvantage of having the
outside sum have a large range as $n$ and $k$ get large.
Nevertheless, we note that for (\ref{HRE}) there can be given a
direct combinatorial proof via a sign-reversing involution so that
it does not require any use of recursions. It is therefore natural
to ask whether one can find similar proofs for our formulas in
sections 3 and 4.

There are several ways in which one could extend our research. For
example, one can study our refined statistics
($\overleftarrow{\Des}_X(\pi)$, $\overleftarrow{\Ris}_X(\pi)$,
$\Lev_X(\pi)$) on the set of all words avoiding a fixed pattern or a
set of patterns (see~\cite{B,BM1,BM2,BM} for definitions of
``patterns in words'' and results on them). More generally, instead
of considering the set of all words, one can consider a subset of it
defined in some way, and then to study the refined statistics on the
subset. Also, instead of considering refined descents, levels, and
rises (patterns of length 2), one can consider patterns of length 3
and more in which the equivalence class of the first letter is
fixed, or, more generally, in which the equivalence classes of more
than one letter (possibly all letters) are fixed. Once such a
pattern (or set of patterns) is given, the questions on avoidance
(or the distribution of occurrences) of the pattern in words over
$[k]$ can be raised.


\end{document}